\RequirePackage{ifpdf}
\ifpdf 
\documentclass[pdftex]{sigma}
\else
\documentclass{sigma}
\fi

\begin{document}

\allowdisplaybreaks

\renewcommand{\thefootnote}{$\star$}

\renewcommand{\PaperNumber}{045}

\FirstPageHeading

\ShortArticleName{The Explicit Construction of Einstein Finsler Metrics
with Non-Constant Flag Curvature}

\ArticleName{The Explicit Construction of Einstein Finsler Metrics\\
with Non-Constant Flag Curvature\footnote{This paper is a
contribution to the Special Issue ``\'Elie Cartan and Dif\/ferential Geometry''. The
full collection is available at
\href{http://www.emis.de/journals/SIGMA/Cartan.html}{http://www.emis.de/journals/SIGMA/Cartan.html}}}

\Author{Enli GUO~$^\dag$, Xiaohuan MO~$^\ddag$ and Xianqiang
ZHANG~$^\S$}

\AuthorNameForHeading{E.~Guo, X.~Mo and X.~Zhang}

\Address{$^\dag$~College of Applied Science, Beijing University of
Technology, Beijing 100022, China}
\EmailD{\href{mailto:guoenli@bjut.edu.cn}{guoenli@bjut.edu.cn}}

\Address{$^\ddag$~Key Laboratory of Pure and Applied Mathematics,
School of Mathematical Sciences,\\
\hphantom{$^\ddag$}~Peking  University, Beijing 100871, China}
\EmailD{\href{mailto:moxh@pku.edu.cn}{moxh@pku.edu.cn}}

\Address{$^\S$~Tianfu College, Southwestern University of Finance and Economics,
Mianyang 621000, China}
\EmailD{\href{mailto:zxq@tf-swufe.net}{zxq@tf-swufe.net}}

\ArticleDates{Received December 08, 2008, in f\/inal form April 09,
2009; Published online April 14, 2009}

\Abstract{By using the Hawking Taub-NUT metric, this note gives an
explicit construction of a 3-parameter family of Einstein Finsler
metrics of non-constant f\/lag curvature in terms of navigation
representation.}

\Keywords{Finsler manifold; Einstein Randers metric; Ricci
curvature}

\Classification{58E20}

\section{Introduction}

The Ricci curvature of a Finsler metric $F$ on a manifold is a
scalar function ${\rm Ric}: TM\rightarrow \mathbb{R}$ with the homogeneity
${\rm Ric}(\lambda y)=\lambda^{2}{\rm Ric}(y)$. See \eqref{(3)} below. A~Finsler
metric $F$ on an $n$-dimensional manifold $M$ is called an
\textit{Einstein metric} if there is a scalar function $K=K(x)$ on
$M$ such that
\[
{\rm Ric}=(n-1)KF^{2}.
\]

Recently, C. Robles studied a special class of Einstein Finsler
metrics, that is, Einstein Randers metrics and obtained the
following interesting result \cite[Proposition 12.9]{ref1}: Let $F$ be a~Randers metric on a $3$-dimensional manifold. Then $F$ is Einstein
if and only if it has constant f\/lag curvature. Note that the f\/lag
curvature in Finsler geometry is a natural extension of the
sectional curvature in Riemannian geometry and it furnishes the
lowest order term for the Jacobi equation that governs the second
variation of geodesics of the Finsler metric. Together with the
classif\/ication theorem of Randers metrics of constant f\/lag
curvature
 due to Bao--Robles--Shen~\cite{ref2}, three dimensional Einstein Randers metrics
  are completely determined.

The Randers metrics were introduced by physicist Randers in 1941
by modifying a Riemannian metric $\alpha:=\sqrt{a_{ij}(x)y^iy^j}$
by a linear term $\beta:=b_i(x)y^i$~\cite{ref8}. By requiring
$\|\beta\|_{\alpha}<1$, we ensure that $\alpha+\beta$ is positive
and strongly convex. The interested reader is referred to~\cite{ref1} for
this result and a thorough treatment of Randers metrics.

  The next problem is to describe four dimensional Einstein
  Randers metrics. This problem turns out to be very dif\/f\/icult. The very f\/irst step
  might be to construct as many examples as possible. In 2002 D.~Bao and C.~Robles
  constructed for the f\/irst time a family of non-constant f\/lag curvature Einstein Randers
  metrics
  on $\mathbb{C}P^{2}$ using Killing vector f\/ields with respect to its Fubini--Study metric~\cite{ref1}.

  The main technique in \cite{ref1} is described as follows. Given a Riemannian
  metric $g={g_{ij}(x)y^{i}y^{j}}$ and a vector f\/ield $W=W^i\frac{\partial}{\partial x^i}$ on a manifold $M$
  with $g(x, W_x)<1$, one can def\/ine a Finsler metric
  $F: TM\rightarrow[0,\infty)$ by
 \begin{gather}\label{(1)}
 F(x,y)=\sqrt g\big(x,y-F(x,y)W_{x}\big).
 \end{gather}
 Solving \eqref{(1)} for $F$, one obtains $F=\alpha+\beta,$ where
 $\alpha=\sqrt{a_{ij}(x)y^{i}y^{j}}$ and $\beta=b_{i}(x)y^{i}$ are
 given by
\begin{gather}\label{(2)}
 a_{ij}=\frac{g_{ij}}{\lambda}+\frac{W_{i}W_{j}}{\lambda^{2}},\qquad b_{j}=-\frac{W_{j}}{\lambda},
 \end{gather}
  where $W_i=g_{ij}W^j$ and $\lambda=1-W_iW^i$.

  In this paper, using \eqref{(1)} we are going to construct a $3$-parameter family of Einstein Randers metrics
  with non-constant f\/lag
  curvature.

  Now let us describe our construction.
Let $(N^3,h)$ be an oriented constant curvature $3$-manifold,
set $M^4=\mathbb{R}\times N^3$, and let $\varphi: M^4\to N^3$ be
the projection onto the second factor. Def\/ine a Riemannian metric
$g$ on $M^4$ by
\[
g=u\varphi^*h+u^{-1}(dt+A)^2,
\]
where $u$ is a positive smooth
function and $A$ a $1$-form on $N^3$. Then $(M^4,g)$ is Einstein
if and only if $u$ and $A$ are related by the monopole equation of
mathematical physics
\[
du=-*dA
\]
and $(N^3,h)$ is f\/lat, in which case $g$ is Ricci-f\/lat, that is,
$g$ has zero Ricci curvature \cite{ref4,ref5}. For $a\geq 0$, def\/ine a
harmonic function on $\mathbb{R}^3\backslash\{0\}$ by
\[
u_a(y)=\frac 14\left(\frac 1{|y|}+a\right).
\]
Then the above construction gives the Hawking Taub-NUT Riemannian
metric $g_a$ $(a>0)$ or the standard metric $g_0$ $(a=0)$. Note
that the metric $g_a$ extends to the whole of $\mathbb{R}^4$; in
fact it is given by the explicit formula~\cite{ref4,ref5,ref7,ref11}
\[
g_{a}=(a|x|^{2}+1)g_{0}-\frac{a(a|x|^{2}+2)}{a|x|^{2}+1}
\left(-x^{2}dx^{1}+x^{1}dx^{2}-x^{4}dx^{3}+x^{3}dx^{4}\right)^{2}.
\]
See \cite{ref5} for discussion of $g_1$. For any $m,n\in \mathbb{R}^+$,
we def\/ine $W_{m,n}\in \Gamma (T\mathbb{R}^4)$ by
\[
W_{m,n}=-mx^2\frac{\partial}{\partial
x^1}+mx^1\frac{\partial}{\partial
x^2}-nx^4\frac{\partial}{\partial
x^3}+nx^3\frac{\partial}{\partial x^4}.
\]

  Let $\Omega:=\{x \in \mathbb{R}^{4}\, |\,f(x)<1\}$, where $f(x)$
  is def\/ined in \eqref{(29)}. We obtain the following result:

 \begin{theorem}\label{theorem1}
  Let $F=\sqrt{a_{ij}(x)y^{i}y^{j}}+b_{i}(x)y^{i}$ be any function in
  $T\Omega\rightarrow [0,\infty)$ on $T\Omega$ which is expressed by \eqref{(2)} in terms of the Hawking
  Taub-NUT metric $g_{a}$ and vector field $W_{m,n}$. Then $F$ has the
  following properties:
\begin{enumerate}\itemsep=0pt
\item[$(i)$] $F$ is a Randers metric;

\item[$(ii)$] $F$ is Einstein with Ricci constant zero;

\item[$(iii)$] $F$ has non-constant flag curvature.
\end{enumerate}
  \end{theorem}

  Our main approach is to show $W_{m,n}$ is the vector f\/ield induced by a one
  parameter isometric group with respect to the Hawking Taub-NUT metric $g_a$.

  These examples show the existence of a large family of global Einstein
  Finsler metrics on $\mathbb{R}^4$ (taking $m=n<a$, then
  $\Omega=\mathbb{R}^4$). This is still far from a complete
  description, of course, but it gives an indication that this family
  is much larger than previously believed.

  \section{Preliminaries}\label{section2}

Let $F$ be a Finsler metric on an $n$-dimensional  manifold $M$.
The second variation of geodesics gives rise to
 a family of endomorphisms
$R_{y}=R^{i}_{\; k}dx^{k}\otimes\frac{\partial}{\partial x^{i}}:
T_{x}M\rightarrow T_{x}M,$ def\/ined by
\[
R_{\; k}^{i}=2\frac{\partial G^{i}}{\partial
x^{k}}-y^{j}\frac{\partial^{2}G^{i}}{\partial x^{j}\partial
y^{k}}+2G^{j}\frac{\partial^{2}G^{i}}{\partial y^{j}\partial
y^{k}}-\frac{\partial G^{i}}{\partial y^{j}}\frac{\partial
G^{j}}{\partial y^{k}},
\]
where $G^{i}$ are the geodesic coef\/f\/icients of $F$ \cite{ref6,ref10}. $F$
is said to be of \emph{constant flag curvature} $K=\lambda,$ if
$R^{i}_{\;k}=\lambda(F^{2}\delta^{i}_{\; k}-FF_{y_{k}}y^{i})$ where
$F_{y^{k}}:=\frac{\partial F}{\partial y^{k}}$ \cite{ref6}. Finsler
metrics of constant f\/lag curvature are the natural extension of
Riemannian metrics of constant sectional curvature.

The Ricci curvature ${\rm Ric}$ is def\/ined to be the trace of $R_{y}$
\begin{gather}\label{(3)}
 {\rm Ric}(y):=R^{k}_{\; k}(x,y).
\end{gather}
${\rm Ric}$ is a well-def\/ined scalar function on $TM\backslash\{0\}.$
$F$ is called an {\em Einstein metric} if there is a~scalar
function $K=K(x)$ on $M$ such that ${\rm Ric}=(n-1)KF^2$.

Consider a Randers metric $F=\alpha+\beta$ on $M$. Let
$g={g_{ij}y^{i}y^{j}}$ be the Riemannian metric and
$W=W^{i}\frac{\partial}{\partial x^{i}}$ be the vector f\/ield on
$M$ such that $F$ is def\/ined by \cite{ref1}. Then D.~Bao and C.~Robles
showed
the following result~\cite{ref1,ref9}:

\begin{proposition}\label{proposition1}
 $F$ is Einstein with Ricci scalar ${\rm Ric}(x):=(n-1)K(x)$ if and only if
\begin{enumerate}\itemsep=0pt
\item[$(i)$] $ W$ satisfies $W_{i|j}+W_{j|i}=-4cg_{ij}$, where
$W_{i}=g_{ij}W^{j}, $ and

\item[$(ii)$] $g$ is an Einstein metric, i.e.\
$\widetilde{\rm Ric}=(n-1)\{K(x)+c^{2}\}g$, where $c={\rm const}$ and
$\widetilde{\rm Ric}$ is the Ricci curvature tensor of $g$, in
particular $K(x)={\rm const}$  if $n\geq 3$.
\end{enumerate}

\end{proposition}

\section{One parameter transformation group}\label{section3}

We rewrite the Hawking Taub-NUT metric $g_{a}$ as
\[
g_{a}=g_{ij}dx^{i}dx^{j},
\]
where
\begin{gather}\label{(4)}
(g_{ij})=G(x)=\left(\begin{array}{cccc}  B-A(x^{2})^{2}& Ax^{1}x^{2}&-Ax^{2}x^{4}& Ax^{2}x^{3}\\
                                       Ax^{1}x^{2}& B-A(x^{1})^{2}& Ax^{1}x^{4}&-Ax^{1}x^{3}\\
                                     -Ax^{2}x^{4}& Ax^{1}x^{4}&B-A(x^{4})^{2}& Ax^{3}x^{4}\\
                                       Ax^{2}x^{3}&-Ax^{1}x^{3}& Ax^{3}x^{4}&B-A(x^{3})^{2}
                                      \end{array}
                                      \right)
                                      \end{gather}
 and
 \begin{gather}\label{(5)}
 B=B(x)=a|x|^{2}+1,\qquad A=A(x)=a\left(1+\frac{1}{B}\right).
 \end{gather}
For any $m, n\in \mathbb{R}^{+},$ we def\/ine $\phi : \mathbb{R}^{4}\to
\mathbb{R}^{4}$ by
\begin{gather}\label{(6)}
\phi_{\theta}=\big(\phi^{1}_{\theta},
\phi^{2}_{\theta},\phi^{3}_{\theta},\phi^{4}_{\theta}\big),
\qquad
\left(\begin{array}{l}  y^1\\
y^2\\
y^3\\
y^4
\end{array}
\right) =\phi_{\theta}^{T} =A_{\theta}
\left(\begin{array}{l}  x^1\\
x^2\\
x^3\\
x^4
\end{array}
\right),
\end{gather}
where $A_{\theta}$ is given by
\begin{gather}\label{(7)}
A_{\theta}=\left(\begin{array}{cc}  A_{\theta,m}& ,0\\
                                      0 & A_{\theta,n}
                                      \end{array}
                                      \right)
\end{gather}
and
\begin{gather}\label{(8)}
A_{\theta,m}=\left(\begin{array}{cc}  \cos(m\theta)&-\sin(m\theta)\\
                                      \sin(m\theta)& \cos(m\theta)
                                      \end{array}
                                      \right),\qquad
A_{\theta,\,n}=\left(\begin{array}{cc}  \cos(n\theta)&-\sin(n\theta)\\
                                      \sin(n\theta)& \cos(n\theta)
                                      \end{array}
                                       \right).
\end{gather}
It is easy to see that
\begin{gather}\label{(9)}
|y|^{2}:=\sum^{4}_{i=1}\left(y^{i}\right)^{2}=\sum^{4}_{i=1}\left(x^{i}\right)^{2}=|x|^{2}.
\end{gather}
Furthermore, $\phi_{\theta}$ is a one-parameter transformation
group.

\section[Killing fields]{Killing f\/ields}\label{section4}

In this section, we explicitly construct a two-parameter family of
Killing f\/ields of  the Hawking Taub-NUT Riemannian metrics. First,
we prove that $\phi_{\theta}$ in \eqref{(6)} is an isometry. It is
equivalent to prove that
\begin{gather}\label{(10)}
A_{\theta}^TG(y)A_{\theta}=G(x).
\end{gather}
Let
\begin{gather}\label{(11)}
H(x):=\left(\begin{array}{cc}  H_1(x)&H_2(x)\\
                                      H_2^T(x)&H_3(x)
                                      \end{array}
                                      \right),
\end{gather}
where
\begin{gather}
H_1(x):=\left(\begin{array}{c}   x^2\\
-x^1
\end{array}
\right)\left(x^2, -x^1\right),
\qquad
H_2(x):=\left(\begin{array}{c}   x^2\\
-x^1
\end{array}
\right)\left(x^4, -x^3\right),\nonumber
\\
\label{(14)}
H_3(x):=\left(\begin{array}{c}   x^4\\
-x^3
\end{array}
\right)\left(x^4, -x^3\right).
\end{gather}
It follows from \eqref{(4)}, \eqref{(11)}
and \eqref{(14)} that
\begin{gather}\label{(15)}
G(x)=B(x)I-A(x)H(x).
\end{gather}
Plugging \eqref{(9)} into \eqref{(5)} yields
\begin{gather}\label{(16)}
A(x)=A(y),\qquad B(x)=B(y).
\end{gather}
Note that $A_{\theta}$ satisf\/ies  that
\begin{gather}\label{(17)}
A_{\theta}^T=A_{\theta}^{-1}.
\end{gather}
By using \eqref{(15)}, \eqref{(16)} and \eqref{(17)}, we obtain that \eqref{(10)} is equivalent to
\begin{gather}\label{(18)}
H(x)=A_{\theta}^TH(y)A_{\theta}.
\end{gather}
 In order to check \eqref{(18)}, from \eqref{(7)}, \eqref{(8)} and \eqref{(11)}, it is enough to check
$3$ matrix equations of order $2\times 2$ as follows:
\begin{gather}
H_1(x)=A_{\theta, m}^TH_1(y)A_{\theta, m},
\qquad
H_2(x)=A_{\theta,m}^TH_2(y)A_{\theta, n},
\qquad
\label{(21)}
H_3(x)=A_{\theta,n}^TH_3(y)A_{\theta, n},
\end{gather}
where $A_{\theta,m}$ and $A_{\theta,n}$ are def\/ined in \eqref{(8)}.
Hence{\samepage
\begin{gather*}
A_{\theta, m}^TH_1(y)A_{\theta, m}=A_{\theta, m}^T\left(\begin{array}{c}   y^2\\
-y^1
\end{array}
\right)\left(y^2, -y^1\right)A_{\theta, m}=\left(\begin{array}{c}   x^2\\
-x^1
\end{array}
\right)\left(x^2, -x^1\right)=H_1(x),
\end{gather*}
thus we obtain the f\/irst equation of \eqref{(21)}, the others are completely analogous.}

For any f\/ixed $p=(p^{1},p^{2},p^{3},p^{4})\in \mathbb{R}^{4},$ we
have
\begin{gather*}
\frac{d\phi_{\theta}^{1}(p)}{d\theta}=-m(p^{1}\sin(m\theta)+p^{2}\cos(m\theta))=-m\phi_{\theta}^{2}(p),
\\
\frac{d\phi_{\theta}^{2}(p)}{d\theta}=m(p^{1}\cos(m\theta)-p^{2}\sin(m\theta))=m\phi_{\theta}^{2}(p).
\end{gather*}
Similarly, we have
\begin{gather*}
\frac{d\phi_{\theta}^{3}(p)}{d\theta}=-n\phi_{\theta}^{3}(p),\qquad
\frac{d\phi_{\theta}^{4}(p)}{d\theta}=n\phi_{\theta}^{4}(p).
\end{gather*}
It
follows that the vector f\/ield induced by $\phi_{\theta}$ is given
by
\begin{gather}
W_{p} = \frac{d}{d\theta}[\phi_{\theta}(p)]\Big|_{\theta=0}
 = \frac{d\phi_{\theta}^{j}(p)}{d\theta}\Big|_{\theta=0}
\frac{\partial}{\partial x^{j}}\Big|_{p}\nonumber\\
\phantom{W_{p}}{}
=-mp^{2}\frac{\partial}{\partial
x^{1}}\Big|_{p}+mp^{1}\frac{\partial}{\partial x^{2}}\Big|_{p}
-np^{4}\frac{\partial}{\partial
x^{3}}\Big|_{p}+np^{3}\frac{\partial}{\partial x^{4}}\Big|_{p}.\label{(22)}
\end{gather}
Note that $\phi_{\theta}$ is an isometry. It follows that $W$ is
of Killing type with respect to the Hawking Taub-NUT
 metric~$g_{a}$.

\section[Construction of Einstein-Finsler metrics]{Construction of Einstein--Finsler metrics}\label{section5}

We rewrite our Killing f\/ield $W$ as
\begin{gather*}
W=W_{m,n}=\sum^{4}_{j=1}W^{j}\frac{\partial}{\partial x^{j}}.
\end{gather*}
Then, from \eqref{(22)}, we have
\begin{gather*}
W^{1}=-mx^{2}, \qquad W^{2}=mx^{1}, \qquad W^{3}=-nx^{4}, \qquad
W^{4}=nx^{3}.
\end{gather*}
Together with \eqref{(4)} we get
\begin{gather*}
W_{j}=g_{ji}W^{i}=\left\{
\begin{array}{ll}
\displaystyle W ^{j}\left (B-\frac{\sigma A}{m}\right),& j=1,2,\vspace{2mm}\\
\displaystyle W^{j}\left(B-\frac{\sigma A}{n}\right),&j=3,4,
\end{array}
\right.
\end{gather*}
where
\begin{gather}\label{(26)}
\sigma:=m\left(x^{1}\right)^{2}+m\left(x^{2}\right)^{2}+n\left(x^{3}\right)^{2}+n\left(x^{4}\right)^{2}
\end{gather}
and
$A$ and $B$ are def\/ined in \eqref{(5)}. It follows that
\begin{gather*}
\lambda=1-|W|^{2},
\end{gather*}
where
\begin{gather}\label{(28)}
|W|^{2}=\sum^{4}_{j=1}W_{j}W^{j}=\left[\left(W^{1}\right)^{2}+\left(W^{2}\right)^{2}\right]\left(B-\frac{\sigma
A}{m}\right)+\left[\left(W^{3}\right)^{2}+\left(W^{4}\right)^{2}\right]\left(B-\frac{\sigma
A}{n}\right).\!\!\!
\end{gather}

We are going to f\/ind the suf\/f\/icient condition producing
Einstein--Finsler metrics in terms of navigation representation,
i.e.\ $|W|<1$.

Put
\[
p:=\max\{m, n\},\qquad q:=\min\{m, n\},
\qquad \mbox{then} \qquad
p-q=|m-n|.
\]
By \eqref{(26)}, we obtain
\[
\frac{\sigma}{m}\geq \frac{q}{p}|x|^{2},\qquad \frac{\sigma}{n}\geq
\frac{q}{p}|x|^{2}.
\]
Thus we have
\begin{gather}
\left[\left(W^{1}\right)^{2}+\left(W^{2}\right)^{2}\right]\left(B-\frac{\sigma
A}{m}\right)+\left[\left(W^{3}\right)^{2}+\left(W^{4}\right)^{2}\right]\left(B-\frac{\sigma
A}{n}\right)\nonumber\\
 \qquad {}\leq
p|x|^{2}\left(1+a|x|^{2}-\frac{aq}{p}|x|^{2}\frac{2+a|x|^{2}}{1+a|x|^{2}}\right)
\nonumber\\
\qquad{}=\frac{|x|^{2}}{1+a|x|^{2}}\left(p+2a|m-n| |x|^{2}+a^{2}|m-n| |x|^{4}\right):=f(x).\label{(29)}
\end{gather}

\begin{proof}[Proof of Theorem~\ref{theorem1}] From \eqref{(28)} and \eqref{(29)}, we obtain that
$f(x)<1$ implies $|W_{m, n}|<1.$ It ensures that
$(g_a, W_{m,n})$ produces Randers metrics in terms of
navigation representation. It is easy to check that
$\sum\limits^{4}_{i=1}\tilde{F}_{x^{i}y^{1}}y^{i}\neq \tilde{F}_{x^{1}},$
where $\tilde{F}:=\sqrt{g_{ij}(x)y^{i}y^{j}}$ and
$g_{ij}=g_a\big(\frac{\partial}{\partial
x^i}, \frac{\partial}{\partial x^j}\big)$. According to~\cite{ref3}, $g_{a}$
is not locally projectively f\/lat. Recall that a Finsler metric
$\bar{F}$ on a $k$-dimensional manifold $M$ is {\em projectively
flat} if every point of $M$ has a neighborhood $U$ that can be
embedded into $\mathbb{R}^k$ in such a way  that it carries the
$\bar{F}$-geodesics in $U$ to straight line segments. In the
Riemannian case, by a theorem of Bonnet--Beltrami, projective
f\/latness is equivalent to having constant sectional curvature.
Hence $g_a$ does not have constant sectional curvature. Note that
$F:=\alpha+\beta$ has constant f\/lag curvature if and only if $g_a$
has constant sectional curvature and $W_{m,n}$ is a homothetic
vector f\/ield~\cite{ref2}. It follows that $F$ does not have constant f\/lag
curvature and we obtain $(iii)$ of Theorem~\ref{theorem1}. Note that the
Hawking Taub-NUT metric $g_{a}$ on $\mathbb{R}^4$ is an Einstein
metric for all $a>0$ \cite{ref7,ref11}. From Section~\ref{section4} we see that
$W_{m,n}$ is of Killing type with respect to $g_a$ for all
$(m,n)\in (\mathbb{R}^+,\mathbb{R}^+)$. Now $(ii)$ of Theorem~\ref{theorem1} is an immediate consequence of Proposition~\ref{proposition1}.
\end{proof}

\pdfbookmark[1]{References}{ref}
\LastPageEnding

\end{document}